\documentclass[12pt, reqno]{amsart}
\usepackage{amssymb, amsmath, amsthm, cancel, hyperref, color, enumerate, verbatim, mathrsfs, mathtools, gensymb, caption, subcaption, graphicx, multicol, setspace, enumitem}

\usepackage[margin=1in]{geometry}

\newtheorem{thm}{Theorem}[section]
\newtheorem{lem}[thm]{Lemma}
\newtheorem{prop}[thm]{Proposition}

\theoremstyle{definition}

\theoremstyle{remark}
\newtheorem*{rem}{Remark}

\newcommand{\CC}{\mathbb{C}}
\newcommand{\ZZ}{\mathbb{Z}}

\newcommand{\SL}{\text{SL}_2(\ZZ)}
\newcommand{\mat}[4]{\begin{bmatrix} #1 & #2 \\ #3 & #4 \end{bmatrix}}

\title{Cusp Forms Without Complex Multiplication as $p$-Adic Limits}
\author{Dalen Dockery}
\begin{document}

\begin{abstract}
In 2016, Ahlgren and Samart used the theory of holomorphic modular forms to obtain lower bounds on $p$-adic valuations related to the Fourier coefficients of three cusp forms. In particular, their work strengthened a previous result of El-Guindy and Ono which expresses a cusp form as a $p$-adic limit of weakly holomorphic modular forms. Subsequently, Hanson and Jameson extended Ahlgren and Samart's result to all one-dimensional cusp form spaces of trivial character and having a normalized form that has complex multiplication. Here we prove analogous $p$-adic limits for several one-dimensional cusp form spaces of trivial character but whose normalized form does not have complex multiplication. 
\end{abstract}

\maketitle

\section{Introduction}
Expressing holomorphic cusp forms as $p$-adic limits of weakly holomorphic modular forms has been of recent interest. In \cite{ElGuindyOno} El-Guindy and Ono considered the cusp form with complex multiplication (CM)
\[g(z) \coloneqq \eta(4z)^2\eta(8z)^2 = \sum\limits_{n=1}^{\infty} a(n)q^n = q-2q^5 -3q^9+6q^{13} + \dots \]
and the weakly holomorphic modular form 
\[F(z) \coloneqq -\frac{\eta(4z)^2\eta(16z)^6}{\eta(32z)^4} = \sum\limits_{n=-1}^{\infty} C(n) q^n = -\frac{1}{q} +2q^3 +q^7-2q^{11}+\cdots.\]
El-Guindy and Ono proved
\begin{equation}\label{eq:EO_limit}
\lim\limits_{m \rightarrow \infty} \frac{F \mid U(p^{2m+1})}{C(p^{2m+1})} = g,
\end{equation}
as a $p$-adic limit, for all primes $p \equiv 3 \pmod{4}$ with $p \nmid C(p).$ Their proof relied on previous work of Guerzhoy, Kent, and Ono (\cite{GKO}) on harmonic Maass forms, studying the $p$-adic coupling of mock modular forms and their shadows.

Ahlgren and Samart strengthened the $p$-adic limit of El-Guindy and Ono, proving the following theorem.

\begin{thm}[Theorem 1.1 of \cite{AhlgrenSamart}]\label{thm:ahlg_sam}
For all $m \geq 0$ and all primes $p \equiv 3 \pmod{4}$, we have
\begin{equation}\label{eq:v(C)}
v_p (C(p^{2m+1})) = m,
\end{equation}
\begin{equation}\label{eq:v(F)}
v_p \left(\frac{F \mid U(p^{2m+1})}{C(p^{2m+1})} - g \right) \geq m+1,
\end{equation}
where $v_p(\cdot)$ denotes the $p$-adic valuation of $\ZZ[[q]].$ 
\end{thm}

\begin{rem}
The case $m=0$ of \eqref{eq:v(C)} proves that $p \nmid C(p)$ for all primes $p \equiv 3 \pmod{4},$ thus removing it as an assumption.
\end{rem}
A key ingredient in their proof is that $g(z)$ lies in a one-dimensional space of cusp forms. Ahlgren and Samart also proved similar results for two other pairs of forms $(F(z),g(z))$. It is noteworthy that Ahlgren and Samart's proof used the theory of classical modular forms and not harmonic Maass forms. 

Hanson and Jameson (\cite{HansonJameson}) furthered the work of Ahlgren and Samart, proving that $p$-adic behavior identical to \eqref{eq:v(C)} and \eqref{eq:v(F)} holds for \textit{all} normalized cusp forms (i.e., having leading coefficient $1$) of one-dimensional cusp form spaces with trivial character. Specifically, Hanson and Jameson proved the following theorem. 

\begin{thm}[Theorem 1 of \cite{HansonJameson}]\label{thm:hanson/jameson}
Let $S_k(N)$ be a one-dimensional space such that the unique normalized cusp form $g(z) = \sum_{n \geq 1} a(n) q^n \in S_k(N)$ has complex multiplication. Then there exists a corresponding weakly holomorphic modular form $F(z) = \sum_{n \geq -1} C(n) q^n \in S_k^\infty$ such that \eqref{eq:v(C)} and \eqref{eq:v(F)} hold for all primes $p$ which are inert in the field of complex multiplication of $g(z).$
\end{thm}

In light of their work, it is natural to investigate those one-dimensional spaces of cusp forms whose normalized form does not have CM. To this end, we now state our main theorem. 

\begin{thm}\label{thm:main}
Let $S_k(N)$ be a one-dimensional space such that $X_0(N)$ has genus $g_N = 0$, and let $g(z) = \sum_{n \geq 1} a(n)q^n \in S_k(N)$ be the unique normalized cusp form. If $g(z)$ does not have complex multiplication, then there is a corresponding weakly holomorphic modular form $F(z) = \sum_{n \geq -1} C(n)q^n$ such that 
\begin{equation}\label{eq:v(C)2}
v_p(C(p^m)) = v_p(C(p))
\end{equation}
\begin{equation}\label{eq:v(F)2}
v_p\left(\frac{F \mid T_k(p^m)}{C(p^m)}-g\right) = (k-1)m - v_p(C(p))
\end{equation}
for all $m \geq 1$ and all primes $p \nmid N$ with $p \nmid a(p)$ and $v_p(C(p)) \leq k-1.$ In particular, if $p \nmid N$ is prime such that $p \nmid a(p)$ and $p \nmid C(p),$ then for all integers $m \geq 1$,
\[v_p(C(p^m)) = 0,\]
\[v_p\left(\frac{F \mid T_k(p^m)}{C(p^m)}-g\right) = (k-1)m.\]
\end{thm}

\begin{rem}
In striking contrast to \eqref{eq:v(C)}, \eqref{eq:v(C)2} shows that $p$-adic valuations of the coefficients $C(p^m)$ are completely determined by that of $C(p).$ Because of this disparity, without CM the functions $g(z)$ are no longer $p$-adic limits of $F(z)$ under applications of the $U(p)$-operator. However, replacing the $U(p)$-operator with the Hecke operator $T_k(p^m)$ yields a $p$-adic limit analogous to \eqref{eq:EO_limit}, namely
\begin{equation}\label{eq:limit}
\lim\limits_{m \rightarrow \infty} \frac{F \mid T_k(p^m)}{C(p^m)} = g.
\end{equation}
\end{rem}

\begin{rem}
Similar $p$-adic limits to \eqref{eq:limit} (involving the Hecke operator $T_k(p^n))$ exist in the context of Weierstrass mock modular forms. Specifically, we refer the interested reader to Theorem 1.3 of \cite{AGOR}.
\end{rem}

After fixing notation and recalling standard definitions and facts about spaces of modular forms in Section \ref{sec:background}, we first identify all possible pairs $(k,N)$ satisfying the hypotheses of Theorem \ref{thm:main}, i.e. such that $S_k(N)$ is a one-dimensional space whose normalized form does not have CM and such that $g_N =  0$. This is done in Section \ref{sec:f/phi}, along with the construction of two families $\phi_n$ and $F_m$ of forms for each weight and level. These functions are pivotal to the proof of Theorem \ref{thm:main}, which we provide in Section \ref{sec:proof}.

\section{Acknowledgements}
The author would like to sincerely thank Marie Jameson for her invaluable advice and guidance, as well as Samuel Wilson for their many fruitful conversations and assistance in editing. Moreover, the author would also like to thank the referee for providing valued suggestions which enhanced this manuscript.

\section{Background}\label{sec:background}
Here we recount several well-known facts about spaces of modular forms. Unless otherwise mentioned, most definitions and facts below can be found in a standard text on modular forms, e.g. \cite{DiamondShurman} or \cite{OnoWeb}. For a function $f : \mathbb{H} \rightarrow \CC$ on the complex upper half-plane and a matrix $\gamma = \mat{a}{b}{c}{d} \in \SL,$ set 
\[f(z) \, \rvert_k \, \gamma \coloneqq (cz+d)^{-k} f\left(\frac{az+b}{cz+d}\right).\]
For $k \in \ZZ,$ we define $M_k^{!}(N)$ to be the $\CC$-vector space of functions $f : \mathbb{H} \rightarrow \CC$ that are holomorphic on $\mathbb{H}$ and meromorphic at the cusps of $\Gamma_0(N)$ and which satisfy $f(z) \mid_k \gamma = f(z)$ for all $\gamma = \mat{a}{b}{c}{d} \in \Gamma_0(N)$. Furthermore, let $M_k^{\infty}(N) \subseteq M_k^{!}(N)$ denote the forms that are holomorphic at each cusp except possibly infinity, and its cuspidal subspace of forms vanishing at each cusp except possibly infinity is denoted $S_k^{\infty} (N).$ Finally, we denote by $M_k(N)$ the forms in $M_k^{\infty}(N)$ which are holomorphic at every cusp, with $S_k(N)$ denoting its cuspidal subspace of forms vanishing at every cusp.  

Each $f \in M_k^!(N)$ has a corresponding Fourier series (at infinity)
\[f(z) = \sum\limits_{n \geq n_0} c(n) q^n,\]
where $q \coloneqq \exp(2\pi i z).$
For a positive integer $m$, the $U$- and $V$-operators are defined on Fourier series by
\[\sum c(n) q^n \, \rvert \, U(m) \coloneqq \sum c(mn) q^n\]
\[\sum c(n) q^n \, \rvert \, V(m) \coloneqq \sum c(n) q^{mn}.\]
Furthermore, we let $T_{k}(p^n)$ be the standard Hecke operators on $M_k^!(N),$ whose action on $q$-series is given by 
\[f(z) \mid T_{k}(p^n) = \sum\limits_{j=0}^{n} p^{(k-1)j} f \mid U(p^{n-j}) \mid V(p^j)\]
for all primes $p \nmid N.$ In particular, for $p \nmid N$ prime,
\[f(z) \mid T_{k} (p) = f \mid U(p) +  p^{k-1} f \mid V(p). \]
It is well-known that 
\[T_{k} (p^n) : M_k^\infty(N) \rightarrow M_k^\infty(N)\]
for all primes $p \nmid N$ and all $n \geq 1.$
Moreover, $T_{k}(p^n)$ restricts to the subspace $S_k^\infty(N),$ i.e. $f(z) \mid T_{k}(n) \in S_k^\infty (N)$ whenever $f(z) \in S_k^\infty(N).$ 

Another crucial operator on spaces of modular forms is the $\Theta$-operator, defined by 
\[\Theta \coloneqq \frac{1}{2\pi i} \frac{d}{dz} = q \frac{d}{dq}.\]
A standard fact about the $\Theta$-operator is that 
\[\Theta^{k-1} : M_{2-k}^\infty(N) \rightarrow S_k^\infty(N).\]

The Dedekind $\eta$-function is given by
\[\eta(z) \coloneqq q^{1/24} \prod\limits_{n=1}^{\infty} (1-q^n).\]
The following criterion allows one to determine if an eta-quotient, a function of the form $\prod\limits_{\delta \mid N} \eta(\delta z)^{r_\delta}$ for finitely many non-zero integers $r_d$, is a modular form. 

\begin{lem}[Theorem 1.64 of \cite{OnoWeb}]\label{lem:eta-quotient1}
Suppose that the eta-quotient $\prod\limits_{\delta \mid N} \eta(\delta z)^{r_\delta}$ satisfies
\[\frac{1}{2}\sum\limits_{\delta \mid N} r_{\delta} = k \in \ZZ, \quad \sum\limits_{\delta \mid N} \delta r_\delta \equiv 0 \pmod{24}, \quad \sum\limits_{\delta \mid N} \frac{N}{\delta} r_\delta \equiv 0 \pmod{24}.\]
Then $\prod\limits_{\delta \mid N} \eta(\delta z)^{r_\delta} \in M_k^!(N,\chi).$ Here $\chi$ is the Dirichlet character given by $\chi(n) = \left(\frac{(-1)^ks}{n}\right),$ where $s = \prod\limits_{\delta \mid N} \delta^{r_{\delta}}.$ 
\end{lem}

Furthermore, the order of vanishing of an eta-quotient at a cusp of $\Gamma_0(N)$ is straightforward to compute. 

\begin{rem}
While Lemma \ref{lem:eta-quotient1} is stated in a slightly general way to treat modular forms with character, all forms we will need turn out to have trivial character, i.e. the Dirichlet character $\chi$ given above is trivial modulo $N$. In this case $\chi = \chi_{\text{triv}}$ is omitted from the notation; that is, $M_k^!(N,\chi_{\text{triv}}) = M_k^!(N).$ 
\end{rem}

\begin{lem}[Theorem 1.65 of \cite{OnoWeb}]\label{lem:eta-quotient2}
The order of vanishing of the eta-quotient $\prod\limits_{\delta \mid N} \eta(\delta z)^{r_\delta} \in M_k^!(N)$ at the cusp $\frac{c}{d}$ of $\Gamma_0(N)$ is given by
\[\frac{N}{24} \sum\limits_{\delta \mid N} \frac{\text{gcd}(d,\delta)^2r_\delta}{\text{gcd}(d,\frac{N}{d})d\delta}.\]
\end{lem}

\section{Preliminary Results}\label{sec:f/phi}
As shown by Hanson and Jameson (\cite{HansonJameson}), if $\dim(S_k(N))=1$ then it must be the case that $N \leq 241.5$; in addition, we note that also $k \leq 10$ if $N > 1$ and $k \leq 26$ if $N=1$. Using Sage to quickly inspect these possibilities, we find that there are 29 spaces for which $\dim(S_k(N))=1$. Table \ref{tab:k,N} below lists the pairs $(k,N)$ of weights and levels of these spaces, respectively, along with the unique normalized cusp form $g(z) \in S_k(N)$ and the genus $g_N$ of the modular curve $X_0(N).$ Here $E_{2k}(z)$ is the normalized Eisenstein series of weight $2k,$ $\Delta(z) \coloneqq \frac{E_4^3(z)-E_6^2(z)}{1728}$ is the modular discriminant of weight 12, and $j(z) \coloneqq \frac{E_4^3(z)}{\Delta(z)}$ is Klein's $j$-invariant.

Henceforth, we refer to pairs $(k,N)$ for which the normalized $g(z) \in S_k(N)$ does not have CM as non-CM pairs. Notice that the spaces of genus $1$ are exactly those of weight $2$, so spaces $S_k(N)$ satisfying the hypotheses of Theorem \ref{thm:main} are equivalent to non-CM pairs with $k \geq 4.$ 

\renewcommand{\arraystretch}{1.5}
\begin{table}
\begin{center}
\begin{tabular}{ |c||c|c|c| }
\hline
$(k,N)$ & normalized $g \in S_k(N)$ & CM? & $g_N$ \\ 
\hline \hline
$(2,11)$ & $\eta(z)^2\eta(11z)^2=q-2q^2-q^3+2q^4+\dots$ & no & 1 \\ \hline
$(2,14)$ & $\eta(z)\eta(2z)\eta(7z)\eta(14z) = q-q^2-2q^3+q^4+\dots$ & no &  1 \\ \hline
$(2,15)$ & $\eta(z)\eta(3z)\eta(5z)\eta(15z) = q-q^2-q^3-q^4+\dots$ & no & 1 \\ \hline
$(2,17)$ & $q-q^2-q^4-2q^5+4q^7+\dots$ & no & 1 \\ \hline
$(2,19)$ & $q-2q^3-2q^4+3q^5-q^7+\dots$ & no & 1 \\ \hline
$(2,20)$ & $\eta(2z)^2\eta(10z)^2 = q-2q^3-q^5+2q^7 + \dots$ & no & 1 \\ \hline
$(2,21)$ & $q-q^2+q^3-q^4-2q^5+\dots$ & no & 1 \\ \hline
$(2,24)$ & $\eta(2z)\eta(4z)\eta(6z)\eta(12z) = q-q^3-2q^5+q^9 + \dots$ & no & $1$ \\ \hline
$(2,27)$ & $\eta(3z)^2\eta(9z)^2 = q-2q^4-q^7+5q^{13}+\dots$ & yes & 1 \\ \hline
$(2,32)$ & $\eta(4z)^2\eta(8z)^2 = q-2q^5-3q^9+6q^{13} + \dots$ & yes & 1 \\ \hline
$(2,36)$ &$\eta(6z)^4 = q-4q^7+2q^{13}+8q^{19}+\dots$ & yes & 1 \\ \hline
$(2,49)$ & $q+q^2-q^4-3q^8-3q^9+\dots$ & yes & 1 \\ \hline
$(4,5)$ & $\eta(z)^4\eta(5z)^4=q-4q^2+2q^3+8q^4+\dots$ & no & 0 \\ \hline
$(4,6)$ & $\eta(z)^2\eta(2z)^2\eta(3z)^2\eta(6z)^2=q-2q^2-3q^3+4q^4+\dots$ & no & 0 \\ \hline
$(4,7)$ & $\frac{\eta(z)^5\eta(7z)^5}{\eta(2z)\eta(14z)}+4\eta(z)^2\eta(2z)^2\eta(7z)^2\eta(14z)^2= q-q^2-2q^3-7q^4 + \dots$ & no & 0 \\ \hline
$(4,8)$ & $\eta(2z)^4\eta(4z)^4=q-4q^3-2q^5+24q^7+\dots $ & no & 0 \\ \hline
$(4,9)$ & $\eta(3z)^8 =q-8q^4+20q^7-70q^{13}+\dots$  & yes & 0 \\ \hline
$(6,3)$ & $\eta(z)^6\eta(3z)^6 = q-6q^2+9q^3+4q^4 + \dots$ & no & 0 \\ \hline
$(6,4)$ & $\eta(2z)^{12} = q-12q^3+54q^5-88q^7 + \dots$ & no & 0 \\ \hline
$(6,5)$ & $q+2q^2-4q^3-28q^4+25q^5+\dots$ & no & 0 \\ \hline
$(8,2)$ & $\eta(z)^8\eta(2z)^8 = q-8q^2+12q^3+64q^4+\dots$ & no & 0 \\ \hline
$(8,3)$ & $\eta(z)^6\eta(3z)^4\left(\eta(z)^3+9\eta(9z)^3\right)^2 = q+6q^2 - 27q^3 -92q^4 + \dots$ & no & 0 \\ \hline
$(10,2)$ & $\frac{\eta(2z)^{28}}{\eta(4z)^8}+16\eta(z)^8\eta(2z)^4\eta(4z)^8 = q+16q^2-156q^3+256q^4+\dots$ & no & 0 \\ \hline
$(12,1)$ & $\Delta(z) =q-24q^2+252q^3-1472q^4+\dots $ & no & 0 \\ \hline
$(16,1)$ & $\Delta(z) E_4(z) =q+216q^2-3348q^3+13888q^4+\dots$ & no & 0 \\ \hline
$(18,1)$ & $\Delta(z) E_6(z) =q-528q^2-4284q^3+147712q^4+\dots$ & no & 0 \\ \hline
$(20,1)$ & $\Delta(z) E_4^2(z) =q+456q^2+50652q^3-316352q^4+\dots$ & no & 0 \\ \hline
$(22,1)$ & $\Delta(z) E_4(z)E_6(z) =q-288q^2-128844q^3-2014208q^4+\dots$ & no & 0 \\ \hline
$(26,1)$ & $\Delta(z) E_4^2(z)E_6(z) =q-48q^2-195804q^3-33552128q^4+\dots$ & no & 0 \\ \hline
\end{tabular}
\caption{\label{tab:k,N}List of $(k,N)$ with $\dim(S_k(N)) = 1.$}
\end{center}
\end{table}
\renewcommand{\arraystretch}{1}

The following proposition establishes, for each non-CM pair $(k,N)$ with $k \geq 4$, the existence of two modular families $F_m$ and $\phi_n$. Once these families are constructed, the proof of Theorem \ref{thm:main} may be carried out in a general way, so as to handle all such cases simultaneously.  

\begin{prop}\label{prop:f_phi}
Let $(k,N)$ be a non-CM pair with $k \geq 4$.
\begin{enumerate}[label=(\alph*)]
\item For all integers $n \geq 2,$ there exists $\phi_n \in M_{2-k}^\infty(N) \cap \ZZ((q))$ of the form 
\[ \phi_n = q^{-n} + \sum\limits_{m=-1}^{\infty} A_n(m) q^m.\]

\item For all integers $m \geq -1,$ there exists a unique $F_m \in S_k^\infty(N) \cap \ZZ((q))$ of the form 
\[F_m = -q^{-m} + \sum\limits_{n=2}^{\infty} C_m(n) q^n.\]

\item Let $F \coloneqq F_1 = \sum\limits_{n=-1}^\infty C(n) q^n$ and $p \nmid N$ be prime. Then
\[ F \mid T_k(p^n) = p^{(k-1)n} F_{p^n} + C(p^n)g\]
for all integers $n \geq 0.$ 
\end{enumerate}

\begin{proof}
\begin{enumerate}[label=(\alph*)]
\item We begin by defining $\phi_2$ explicitly for each pair $(k,N)$, as well as a modular function $L \in M_0^\infty(N)$ with a simple pole at infinity. We then define $\phi_n$ for $n \geq 3$ inductively by taking appropriate linear combinations of $\phi_{n-1}L, \phi_{n-1},\dots,\phi_2.$ For example, for the pair $(8,2)$ we have
\begin{align*}
\phi_3 &= \phi_2L +104\phi_2 = q^{-3}-2980q^{-1}+146432 - 3896490q+\dots \\
\phi_4 &= \phi_3L +24\phi_3 + 2704\phi_2 = q^{-4} - 71936q^{-1} +3793608 - 12159360q + \dots,
\end{align*}
and so on. 
Table \ref{tab:phi's_F's_kneq2} gives an explicit description of all required $L$ and $\phi_2$. Most of these functions are eta-quotients, in which case their modular behavior is readily verified by Lemmas \ref{lem:eta-quotient1} and \ref{lem:eta-quotient2}. For the remaining pairs $(k,N)$ the functions $\phi_2 \in M_{2-k}^\infty(N)$ are given as products of eta-quotients and sums involving the unique normalized cusp form $g \in S_k(N)$ and the Eisenstein series $E_2(z),E_2(Nz),E_4(z)$ and $E_4(Nz)$; as such, the modularity of these functions is also immediate.

\item Similarly, let $F_{-1} \coloneqq -g$ and inductively define $F_m$ for $m \geq 0$ as linear combinations of $F_{n-1}L, F_{n-1},\dots,F_{-1}.$ For example, in the case of the pair $(8,2)$
\begin{align*}
F_0 &= F_{-1}L+32F_{-1} = -1 -224q^2 + 4096q^3 -31200q^4 + \dots \\
F_1 &= F_0 L +24F_0 - 500F_{-1} = -q^{-1} + 2144q^2-98226q^3+1817856q^4+\dots,
\end{align*}
and so on.

To show uniqueness, suppose two such $F_m$ exist, say $F_m$ and $F_m'$. Then $F_m-F_m'$ must be of the form $O(q^2) \in S_k(N),$ but $g = q +O(q^2)$ spans $S_k(N),$ and thus $F_m-F_m'=0.$ 

\item The proof is identical to that of Proposition 3(c) of \cite{HansonJameson}.
\end{enumerate}
\end{proof}
\end{prop}

\renewcommand{\arraystretch}{2}
\begin{table}
\begin{center}
\begin{tabular}{ |c||c|c| }
\hline
$(k,N)$ & $L \in M_0^\infty(N)$ & $\phi_2 \in M_{2-k}^\infty(N) $ \\ 
\hline \hline
$(12,1)$ & $j(z) = q^{-1}+744+196884q +\dots$ & $\frac{E_4^2(z)E_6(z)}{\Delta(z)^2} = q^{-2}+24q^{-1}-196560+\dots$ \\ \hline
$(16,1)$ & $j(z) = q^{-1}+744+196884q +\dots$ & $\frac{E_4(z)E_6(z)}{\Delta(z)^2}=q^{-2}-216q^{-1}-146880+\dots $ \\ \hline
$(18,1)$ & $j(z) = q^{-1}+744+196884q +\dots$ & $\frac{E_4^2(z)}{\Delta(z)^2} =q^{-2}+528q^{-1}+86184+\dots $ \\ \hline
$(20,1)$ & $j(z) = q^{-1}+744+196884q +\dots$ & $\frac{E_6(z)}{\Delta(z)^2} = q^{-2}-456q^{-1}-39600+\dots$ \\ \hline
$(22,1)$ & $j(z) = q^{-1}+744+196884q +\dots$ & $\frac{E_4(z)}{\Delta(z)^2} = q^{-2}+288q^{-1}+14904+\dots$ \\ \hline
$(26,1)$ & $j(z) = q^{-1}+744+196884q +\dots$ & $\frac{1}{\Delta(z)^2} = q^{-2}+48q^{-1}+1224+\dots$ \\ \hline
$(8,2)$ & $\frac{\eta(z)^{24}}{\eta(2z)^{24}} = q^{-1}-24+276q+\dots$ & $\frac{\eta(z)^{80}}{\eta(2z)^{64}} = q^{-2}-80q^{-1}+3144-80960q+\dots$ \\ \hline
$(10,2)$ & $\frac{\eta(z)^{24}}{\eta(2z)^{24}} = q^{-1}-24+276q+\dots$ & $\frac{\eta(z)^{16}}{\eta(2z)^{32}} = q^{-2}-16q^{-1}+136+\dots$ \\ \hline
$(6,3)$ & $\frac{\eta(z)^{12}}{\eta(3z)^{12}} =  q^{-1}-12+54q+\dots$ & $\frac{\eta(z)^6(3 E_2(3z) - E_2(z))}{2\eta(3z)^{18}}=q^{-2}+6q^{-1}-27-68q+\dots$ \\ \hline
$(8,3)$ & $\frac{\eta(z)^{12}}{\eta(3z)^{12}} =  q^{-1}-12+54q+\dots$ & $\frac{\eta(z)^6}{\eta(3z)^{18}} = q^{-2}-6q^{-1}+q+28q+\dots$ \\ \hline
$(6,4)$ & $\frac{\eta(z)^8}{\eta(4z)^8} = q^{-1}-8+20q-62q^3+\dots$ & $\frac{\eta(2z)^8}{\eta(4z)^{16}} = q^{-2}-8+36q^2-128q^4+\dots$ \\ \hline
$(4,5)$ & $\frac{\eta(z)^6}{\eta(5z)^6}=q^{-1}-6+9q+10q^2+\dots$ & $\frac{\eta(z)^2(5E_2(5z)-E_2(z))}{4\eta(5z)^{10}}= q^{-2}+4q^{-1}+5-16q+\dots$ \\ \hline
$(6,5)$ & $\frac{\eta(z)^6}{\eta(5z)^6}=q^{-1}-6+9q+10q^2+\dots$ & $\frac{\eta(z)^2}{\eta(5z)^{10}} = q^{-2}-2q^{-1}-1+2q+\dots$ \\ \hline
$(4,6)$ & $\frac{\eta(z)^5\eta(3z)}{\eta(2z)\eta(6z)^5} = q^{-1}-5+6q+\dots$ & $\frac{\eta(2z)^4\eta(3z)^6}{\eta(z)^2\eta(6z)^{12}} = q^{-2}+2q^{-1}+1-4q+\dots$ \\ \hline
$(4,7)$ & $\frac{\eta(z)^4}{\eta(7z)^4} = q^{-1}-4+2q+8q^2+\dots$ & $\frac{\eta(z)^2}{10\eta(7z)^{14}}\left(\frac{E_4(z)-E_4(7z)}{240}-g\right)=q^{-2}+q^{-1}+\dots$  \\ \hline
$(4,8)$ & $\frac{\eta(4z)^{12}}{\eta(2z)^4\eta(8z)^8} = q^{-1}+4q+2q^3+\dots$ & $\frac{\eta(4z)^4}{\eta(8z)^8} = q^{-2}-4q^2+10q^6-24q^{10}+\dots$  \\ \hline
\end{tabular}
\caption{\label{tab:phi's_F's_kneq2}Defining $L$ and $\phi_2$ for $(k,N), k \geq 4.$}
\end{center}
\end{table}
\renewcommand{\arraystretch}{1}

For non-CM pairs with $k \geq 4$, the coefficients of $F_m$ and $\phi_n$ satisfy a property known as Zagier duality. Along with the uniqueness of the functions $F_m$ established in Proposition \ref{prop:f_phi}(b), the duality of these coefficients immediately implies that the functions $\phi_n$ are also uniquely defined in such cases. 

\begin{lem}\label{lem:Zagier_duality}
Let $(k,N)$ be a non-CM pair with $k \geq 4$, and let $\phi_n = q^{-n} + \sum_{m \geq -1} A_n(m) q^m \in M_{2-k}^\infty(N)$ and $F_m = -q^{-m} + \sum_{n \geq 2} C_m(n) q^n \in S_k^\infty(N)$ be as above. Then for all $n \geq 2$ and all $m \geq -1,$ we have 
\[C_m(n) = A_n(m).\]
\begin{proof} Zagier duality for the canonical bases of $S_k^\infty(N)$ and $M_{2-k}^\infty(N)$ has been established for all $N$ of genus $0$. Specifically, we refer to a collection of works by Jenkins and collaborators, such as \cite{JenkinsThornton}, \cite{HaddockJenkins},\cite{IbaJenkinsWarnick}, \cite{JenkinsThornton2}, and \cite{JenkinsMolnar}.
\end{proof}
\end{lem}

\section{Proof of Theorem \ref{thm:main}}\label{sec:proof} 
We are now equipped to prove our main theorem. Fix a non-CM pair $(k,N)$ with $k \geq 4$, and let $p \nmid N$ be a prime such that $p \nmid a(p)$ and $v_p(C(p)) \leq k-1.$ We first show that \eqref{eq:v(C)2} holds, from which \eqref{eq:v(F)2} will quickly follow. 

From Proposition \ref{prop:f_phi}(a), $\phi_p$ is of the form
\[\phi_p = q^{-p} +
\sum\limits_{n=-1}^{\infty} A_p(n) q^n.\]
Since $k$ is even,
\begin{align*}
\Theta^{k-1} (\phi_p) &= -p^{k-1}q^{-p} -A_{p}(-1)q^{-1} + A_{p}(1)q +\sum\limits_{n=2}^\infty n^{k-1} A_p(n)q^n \\ 
&= -p^{k-1}q^{-p} - C_{-1}(p)q^{-1} + C_1(p) q + \sum\limits_{n=2}^\infty n^{k-1} A_p(n)q^n \\
&= -p^{k-1}q^{-p} + a(p) q^{-1} + C(p) q + O(q^2) ,
\end{align*}
by the Zagier duality established in Lemma \ref{lem:Zagier_duality}. From Proposition \ref{prop:f_phi}(c), we also have 
\begin{align*}
F \mid T_k(p) &= p^{k-1}F_p + C(p)g \\
&= -p^{k-1}q^{-p}+C(p)q+O(q^2).
\end{align*}
Hence
\[F \mid T_k(p) - \Theta^{k-1}(\phi_p) = -a(p)q^{-1} + O(q^2) \in S_k^\infty(\Gamma_0(N)),\]
so by the uniqueness established in Proposition \ref{prop:f_phi}(b), this difference must be 
\[a(p)F = a(p)F_1 = -a(p)q^{-1}+O(q^2).\]
Set $a \coloneqq a(p)$ for ease of notation and rearrange so that
\[ F \mid U(p) = \Theta^{k-1}(\phi_p) +aF -p^{k-1} F \mid V(p).\]
Since $F \mid V(p) \mid U(p) = F,$ we see
\begin{align*} 
F \mid U(p^2) &= F \mid U(p) \mid U(p) \\
&= \Theta^{k-1}(\phi_p) \mid U(p) +aF\mid U(p) -p^{k-1} F \mid V(p) \mid U(p) \\
&= \Theta^{k-1}(\phi_p) \mid U(p) + a \left[\Theta^{k-1}(\phi_p) +aF -p^{k-1} F \mid V(p)\right] -p^{k-1} F \\
&= \Theta^{k-1}(\phi_p) \mid U(p) + a \Theta^{k-1}(\phi_p) + (a^2-p^{k-1})F - ap^{k-1}F \mid V(p)
\end{align*}
Continuing to apply $U(p)$ this way, we write
\begin{equation}\label{eq:applying_U}
F \mid U(p^m) = \left[\sum\limits_{i=0}^{m-1} \alpha_{i,m} \Theta^{k-1}(\phi_p) \mid U(p^i)\right] + \beta_m F + \gamma_m F \mid V(p).
\end{equation}
for coefficients $\alpha_{i,m},\beta_m,\gamma_m.$ We now determine $\alpha_{i,m},\beta_m,$ and $\gamma_m,$ as \eqref{eq:v(C)2} will follow by comparing the coefficients of $q$ in \eqref{eq:applying_U}. Observe
\begin{align*} 
F \mid U(p^{m+1}) &= F \mid U(p^m) \mid U(p) = \left[\sum\limits_{i=0}^{m-1} \alpha_{i,m} \Theta^{k-1}(\phi_p) \mid U(p^{i+1}) \right] + \beta_m F \mid U(p) + \gamma_m F \\
&= \left[\sum\limits_{i=1}^{m} \alpha_{i-1,m} \Theta^{k-1}(\phi_p) \mid U(p^{i}) \right] + \beta_m \left(\Theta^{k-1}(\phi_p) +a F - p^{k-1} F \mid V(p) \right) + \gamma_m F \\
&= \beta_m \Theta^{k-1}(\phi_p) + \left[\sum\limits_{i=1}^{m} \alpha_{i-1,m} \Theta^{k-1}(\phi_p) \mid U(p^{i}) \right] + \left(a \beta_m + \gamma_m\right) F - p^{k-1}\beta_m F \mid V(p).
\end{align*}
As such, the coefficients $\alpha_{i,m},\beta_m,$ and $\gamma_m$ satisfy the recurrences
\begin{align}
\alpha_{i,m+1} &= \begin{cases} \beta_m, & i=0 \\ \alpha_{i-1,m}, & i \geq 1 \end{cases} \label{eq:alphas} \\
\beta_{m+1} &= a \beta_m + \gamma_m \label{eq:betas} \\
\gamma_{m+1} &= -p^{k-1}\beta_m. \label{eq:gammas}
\end{align} 
Notice that applying \eqref{eq:alphas} inductively yields $\alpha_{i,m+1} = \beta_{m-i},$ and by substituting \eqref{eq:gammas} into \eqref{eq:betas} we find
\[\beta_{m+1} = a \beta_m - p^{k-1} \beta_{m-1}.\]
This second-order recurrence can be solved to yield the explicit formula
\[\beta_m =  \frac{(a+\sqrt{a^2-4p^{k-1}})^{m+1} -(a-\sqrt{a^2-4p^{k-1}})^{m+1}}{2^{m+1}\sqrt{a^2-4p^{k-1}}}.\] 
From the Binomial Theorem, 
\begin{align*} 
\beta_m &= \frac{1}{2^{m+1}\sqrt{a^2-4p^{k-1}}} \sum\limits_{i=0}^{m+1} {m+1 \choose i} a^{m+1-i} \left(\sqrt{a^2-4p^{k-1}}^i - (-1)^{i} \sqrt{a^2-4p^{k-1}}^i\right) \\ 
&= \frac{1}{2^m \sqrt{a^2-4p^{k-1}}} \sum\limits_{\substack{i=0 \\ i \text{ odd}}}^{m+1} {m+1 \choose i} a^{m+1-i} \sqrt{a^2-4p^{k-1}}^i \\
&= \frac{1}{2^m} \sum\limits_{\substack{i=0 \\ i \text{ even}}}^{m} {m+1 \choose i+1} a^{m-i} (a^2-4p^{k-1})^{i/2} \\
&= \frac{1}{2^m} \sum\limits_{i=0}^{\lfloor \frac{m}{2} \rfloor} {m+1 \choose 2i+1} a^{m-2i} (a^2-4p^{k-1})^i \\
&= \frac{1}{2^m} \sum\limits_{i=0}^{\lfloor \frac{m}{2} \rfloor} {m+1 \choose 2i+1} a^{m-2i} \sum\limits_{j=0}^{i} {i \choose j} (a^2)^{i-j} (-4p^{k-1})^j \\
&= \frac{1}{2^m} \sum\limits_{j=0}^{\lfloor \frac{m}{2} \rfloor} \sum\limits_{i=j}^{\lfloor \frac{m}{2} \rfloor} {m+1 \choose 2i+1} {i \choose j} (-4)^j a^{m-2j} p^{(k-1)j}.
\end{align*}
As such, we may write 
\begin{equation}\label{eq:newbeta}
\beta_m = \sum\limits_{j=0}^{\lfloor \frac{m}{2} \rfloor} c_{m,j} p^{(k-1)j},
\end{equation}
where 
\begin{equation}\label{eq:c}
c_{m,j} \coloneqq \frac{(-4)^ja^{m-2j}}{2^m} \sum\limits_{i=j}^{\lfloor \frac{m}{2} \rfloor} {m+1 \choose 2i+1} {i \choose j}.
\end{equation}

Now we compare the coefficients of $q$ in \eqref{eq:applying_U} to obtain 
\begin{align}
C(p^m) &= \sum\limits_{i=0}^{m-1} \alpha_{i,m} p^{(k-1)i} A_p(p^i) = \sum\limits_{i=0}^{m-1} \beta_{m-i-1} p^{(k-1)i} A_p(p^i) = \sum\limits_{i=0}^{m-1} \beta_i p^{(k-1)(m-i-1)} A_p(p^{m-i-1}) \nonumber \\
&= \sum\limits_{i=0}^{m-1} \sum\limits_{j=0}^{\lfloor \frac{i}{2} \rfloor} c_{i,j}  p^{(k-1)(m+j-i-1)} A_p(p^{m-i-1}) \label{eq:c(p^m)} 
\end{align}
thanks to \eqref{eq:newbeta}.

If $m=1,$ we have nothing to prove, so assume $m  \geq 2.$ Notice that $m+j-i-1 
\geq 1$ for all $0 \leq i \leq m-1$ and $0 \leq j \leq \lfloor \frac{i}{2} \rfloor$, except for when $i=m-1$ and $j=0.$ Hence, all terms in the double sum of \eqref{eq:c(p^m)} are divisible by at least $p^{k-1},$ other than $c_{m-1,0}A_p(1).$ But by \eqref{eq:c}
\[c_{m-1,0} A_p(1)= \frac{a^{m-1}}{2^{m-1}} \sum\limits_{i=0}^{\lfloor \frac{m-1}{2} \rfloor} {m \choose 2i+1}A_p(1) = \frac{a^{m-1}}{2^{m-1}} \cdot 2^{m-1} A_p(1) = a^{m-1}C(p),\] 
using Lemma \ref{lem:Zagier_duality}.  Now since $p \nmid a$ by assumption, 
\[v_p(c_{m-1,0}A_p(1)) = v_p(a^{m-1} C(p)) = v_p(C(p)) \leq k-1.\]
Therefore, $v_p(C(p^m)) = v_p(C(p)),$ which completes the proof of $\eqref{eq:v(C)2}$. 

To establish \eqref{eq:v(F)2}, note that from Proposition \ref{prop:f_phi} 
\[\frac{F \mid T_k(p^m)}{C(p^m)} - g = \frac{p^{(k-1)m} F_{p^m}}{C(p^m)}.\]
Thus clearly 
\[v_p \left(\frac{F \mid T_k(p^m)}{C(p^m)} - g\right) = v_p\left(\frac{p^{(k-1)m}F_{p^m}}{C(p^m)}\right)= (k-1)m-v_p(C(p^m)) = (k-1)m - v_p(C(p)),\]
thanks to \eqref{eq:v(C)2}. Note that equality holds since $p \nmid F_{p^m},$ as the coefficient of $q^{p^m}$ in $F_{p^m}$ is $-1$.

\section{Concluding Remarks}
By coupling Theorem \ref{thm:main} with Hanson and Jameson's main result (Theorem 1 of \cite{HansonJameson}), $p$-adic limits have been established for a large portion of cusp forms in one-dimensional spaces $S_k(N)$. More precisely, only $p$-adic limits for non-CM forms of weight 2 remain to be shown. While such limits appear to hold computationally, our present approach is not equipped to handle these cases. The main difficulty in treating weight 2 forms stems from the potential lack of Zagier duality for levels $N$ with $g_N = 1$. In these cases, for a prime $p$ one has to construct $\phi_p$ of the form 
\[\phi_p = q^{-p} + a(p)q^{-1} + C + C(p)q + O(q^2)\]
for some constant $C$, rather than leveraging duality to guarantee that $A_p(-1) = C_{-1}(p)=a(p)$ and $A_p(1) = C_1(p) = C(p)$. However, once these $\phi_n$ are constructed, the proof technique in Section \ref{sec:proof} will immediately apply. 

One could also investigate cusp forms without CM lying in spaces of dimension larger than 1. For instance, in this direction Dicks (\cite{Dicks}) proved $p$-adic limits for all weight 2 newforms in one-dimensional cuspidal subspaces that are expressible as eta-quotients. It is plausible that Dicks' methods could be adapted to handle newforms without CM, perhaps even in the few remaining $\dim(S_k(N))=1$ cases.
\bibliographystyle{alpha}
\bibliography{refs}
\end{document}